\documentclass[12pt]{article}
\usepackage{amssymb}

\usepackage{amsmath}
\usepackage{tikz}
\usepackage{amscd}
\usepackage{latexsym}
\usepackage{mathrsfs}
\usepackage{enumitem}
\usepackage{lipsum}
\usepackage[noadjust]{cite}

\newtheorem{Theorem}{Theorem}[section]

\newtheorem{Lemma}[Theorem]{Lemma}

\newtheorem{Proposition}[Theorem]{Proposition}
\newtheorem{Corollary}[Theorem]{Corollary}

\newcommand{\qed}{\hfill $\Box$}
\newenvironment{Proof}{ \emph{Proof.}}{\qed}
\def\bc{\begin{center}}
\def\ec{\end{center}}

\numberwithin{equation}{section}

\begin{document}
\title{{\bf A group and the completion of its coset semigroup}\footnote{This paper is supported by National Natural Science Foundation of China (11971383, 12201495), and Shaanxi Fundamental Science Research Project for Mathematics and Physics (22JSY023).}}
\author{{\bf Xian-zhong Zhao}, {\bf Zi-dong Gao}, {\bf Dong-lin Lei}
  \footnote{Corresponding author. E-mail: dllei@nwu.edu.cn}  \\
   {\small School of Mathematics, Northwest University}\\
   {\small Xi'an, Shaanxi, 710127, P.R. China}\\}
\date{}
\maketitle
\vskip -4pt
\baselineskip 16pt

{\small
\noindent\textbf{ABSTRACT:} Let ${\cal K}_1(G)$ denote the inverse subsemigroup of ${\cal K}(G)$ consisting of all right cosets of all non-trivial subgroups of $G$.
This paper concentrates on the study of the group $\Sigma({\cal K}_1(G))$  of all units of the completion of ${\cal K}_1(G)$.
The characterizations and the representations of $\Sigma({\cal K}_1(G))$ are given when $G$ is a periodic group whose minimal subgroups permute with each other.
Based on these, for such groups $G$ except some special $p$-groups, it is shown that $G$ and its coset semigroup ${\cal K}_1(G)$ are uniquely determined by each other, up to isomorphism. This extends the related results obtained by Schein in 1973.

 }
\section{Introduction}
Let $G$ be a group and ${\cal K}(G)$ denote the set of all right cosets of all subgroups of $G$. Then ${\cal K}(G)$ equipped with the binary operation
$$Ha* Kb = \langle H, K^{a^{-1}}\rangle ab$$
becomes an inverse monoid (see \cite{schein1966} or \cite{mcali2}).
An inverse subsemigroup of ${\cal K}(G)$ is said to be a \emph{coset semigroup} of $G$.
The coset semigroups of a group  play an impotent role in the theory of inverse semigroups.
They have attracted the attention and investigation of many scholars such as,
Schein, McAlister, East, Ara\'{u}jo and Shum et al. (see, \cite{schein1966, schei4, mcali2, james, james2, shum, arau}).
It is easy to see that  ${\cal K}_1 (G)\cong {\cal K}_1 (\widetilde{G})$ if $ G\cong \widetilde{G}$ for any two groups $G$ and $\widetilde{G}$, where ${\cal K}_1(G)$ (resp. ${\cal K}_1 (\widetilde{G})$) denotes the inverse subsemigroup of ${\cal K}(G)$ (resp. ${\cal K} (\widetilde{G})$) consisting of all right cosets of all non-trivial subgroups of $G$ (resp. $\widetilde{G}$),
but not vice versa.
Schein \cite{schei2} tells us that an inverse semigroup and its completion are uniquely determined by each other, up to isomorphism. That is to say, for any inverse semigroups $S$ and $T$, $S\cong T$ if and only if $C(S)\cong C(T)$, where $C(S)$ and $C(T)$ are the completion of $S$ and $T$ respectively. 
In this paper we shall devote to the investigation of the groups $G$ for which ${\cal K}_1 (G)\cong {\cal K}_1 (\widetilde{G})$ implies that $ G\cong \widetilde{G}$ for any group $\widetilde{G}$. As a consequence, we obtain that for such groups, the group $G$, the coset semigroup ${\cal K}_1(G)$ of $G$ and the completion $C({\cal K}_1(G))$ of ${\cal K}_1(G)$ are uniquely determined by each other, up to isomorphism. It not only extends the related result obtained by Schein in \cite{schei2}, but also  build a bridge between the investigations of such groups and of their coset semigroups.

The following concepts, notions and results on the completion of an inverse semigroup due to Schein are needed for us.
Following Schein \cite{schei2}, we call two elements $x$ and $y$ of an inverse semigroup $S$ to be \emph{compatible}, written $x\thicksim y$, if both $xy^{-1}$ and $x^{-1}y$ are idempotents. A subset $A$ of $S$ is said to be \emph{compatible} if the elements of $A$ are pairwise compatible. We also know that $S$ equipped with  the natural partial order becomes  a partial order set. A compatible subset $A$ of $S$ is said to be
\emph{permissible} if it is an order ideal. The set of all permissible subsets of $S$  is denoted by $C(S)$, called the \emph{completion} of $S$.
It is well-known that $C(S)$ is a complete and infinitely distributive inverse monoid under the usual multiplication of subsets and that the natural partial order on $C(S)$  coincide with subset inclusion.
An  idempotent of $C(S)$ is just an order ideal contained in  $E(S)$. The identity of $C(S)$ is $E(S)$. For any $A\in C(S)$,  the inverse  $A^{-1}$ of $A$ is equal to $\{a^{-1} \mid a \in A\}$.
Also, $S$ can be embedded in $C(S)$ by means of the embedding $\iota: S\to C(S)$ defined by $\iota(s)=[s]$, where $[s]$ is the order ideal generated by $s$. The group of units of $C(S)$ is denoted by $\Sigma(S)$ (see Schein \cite{schei2} or Lawson \cite{lawson} for details).
There is a series of papers in the literature devoted to the study of the completions of inverse semigroups and their applications.
Schein \cite{schei2} describes the translational hulls and ideal extensions of an inverse semigroup via its completion.
McAlister \cite{mcali5} and Lawson \cite{lawson1} introduce and study the $E$-unitary inverse semigroups and the almost factorisable inverse semigroups $S$ via the semigroup $C(S)$ and the group $\Sigma(S)$.
Also, the semigroup $C(S)$ has numerous applications in the researches of many other related topics, such as $S$-systems, expansions of inverse semigroups and pseudogroups (see, \cite{schei3, shoji,lawson2, lawson3, cast}).

Noticing that the mapping $\varphi: G \rightarrow \Sigma({\cal K}_1(G))$ defined by $\varphi(a)=\{Ha \mid 1 \neq H \leq G\}$ is a group homomorphism whose kernel is the intersection of all non-trivial subgroups of $G$, we can see that a group $G$ is embedded to the unit group $\Sigma({\cal K}_1(G))$ of $C({\cal K}_1(G))$ if the intersection of all non-trivial subgroups of $G$ is trivial. This motivates us to concentrate on the study of the unit group $\Sigma({\cal K}_1(G))$.

Unless otherwise specified in this article, $G$ always denotes a non-trivial periodic group and  $\Omega_G = \{H_i\}_{i\in I}$ the set of all minimal subgroups of $G$, and writes
$${\cal A}_G=\{\{H_ia_i\}_{i\in I}\mid (\forall \,j, k\in I)\,H_ja_j\thicksim H_ka_k \text{ and } |\{H_ia_i\}_{i\in I}\cap L_{\scriptscriptstyle {H_j}}|=1\}.$$
Then ${\cal A}_G$ equipped with the multiplication $A\cdot B=\overline {[A][B]}$ becomes a group, and the mapping $\varphi$ from ${\cal A}_G$ to $ \Sigma({\cal K}_1(G))$  defined by $\varphi (A) = [A]$ is an isomorphism  (see, Lemma \ref{tonggou}).
Also, it is proved that
$$A\cdot B=\{H_ia_ib_{i^{\tau(A)}}\}_{i\in I}$$
for any $A=\{H_ia_i\}_{i\in I}$, $B=\{H_ib_i\}_{i\in I}\in {\cal A}_G$, where $\tau$ is a group homomorphism from ${\cal A}_G$ to the symmetric group $Sym(I)$.

When $H_i H_j = H_j H_i$ for any $H_i,\, H_j \in \Omega_G $, for any given $H \in \Omega_G $, we
introduce a simple graph $\Gamma(H)$ corresponding to $H$ and prove that
$$|\Sigma({\cal K}_1(G))|=|{\cal A}_G|=|G||H|^{|J|-1},$$
where $|J|$  denotes the cardinality of the set of all connected components of $\Gamma(H)$.
Also, we give the value of $|J|$ by means of the related results in group theory  and the theory of vector spaces.

Based on these, for a periodic group $G$ whose the minimal subgroups permute with each other, we shall show that $G$ is isomorphic to $\Sigma({\cal K}_1(G))$ except some special  $p$-groups.
That is to say, for such  groups, $G$ and $\Sigma({\cal K}_1(G))$ are uniquely determined by each other, up to isomorphism.
Also, Schein tells us  in \cite{schei2} that an inverse semigroup and its completion are uniquely determined by each other, up to isomorphism. That is to say, for any inverse semigroups $S$ and $T$, $S$ and $T$ are isomorphic if and only if $C(S)$ and $C(T)$ are isomorphic.
Thus we can see that for the periodic groups $G$ whose the minimal subgroups permute with each other except some special $p$-groups, $G$ and its coset semigroup ${\cal K}_1(G)$ are uniquely determined by each other, up to isomorphism. This extends the related result obtained by Schein in \cite{schei2}.

For other unexplained notation and terminology used in this paper, the reader is referred to \cite{lawson,petrich,robinson}.

\section{Preliminary}

In this section, the order filters and  the compatible relations in  the coset semigroups ${\cal K}_1(G)$ are studied.
Some  properties of ${\cal K}_1(G)$ are given. It is easy to see that the following elementary properties of ${\cal K}(G)$ holds also in ${\cal K}_1(G)$.

\begin{Lemma}\cite[Lemma 1.1]{mcali2}\label{kgjiben}
\rm Let $G$ be a group and let $Ha, Kb\in{\cal K}(G)$. Then
\begin{itemize}
\item[$(i)$] $(Ha)^{-1}=a^{-1}H=H^{a}a^{-1}$; thus $(Ha)*(Ha)^{-1}=H,$ $(Ha)^{-1}*(Ha)=H^{a}$; further $Ha=Kb$ implies $H=K,\,ab^{-1}\in H$;
\item[$(ii)$] $Ha{\,\cal R\,}Kb$ $[Ha{\,\cal L\,}Kb]$ if and only if $H = K$ [$H^a=K^b$];

$Ha{\,\cal D\,}Kb$ if and only if $H$ is conjugate to $K$;

$Ha\leq_{\cal J} Kb$ if and only if $H$ is conjugate to a subgroup of $K$;

\item[$(iii)$] the idempotents of ${\cal K}(G)$ are precisely the subgroups of $G$; the central  idempotents are the normal subgroups of $G$;

\item[$(iv)$] if $H$ is an idempotent of ${\cal K}(G)$ (i.e., $H$ is a subgroup of $G$), then the $\cal H$-class of ${\cal K}(G)$
containing $H$ is isomorphic to $N_G(H)/H$;

\item[$(v)$] the natural partial order on ${\cal K}(G)$ is the inverse of inclusion.
\end{itemize}
\end{Lemma}

In the paper we write $\leq_n$ to denote the natural partial order on ${\cal K}(G)$.
It is well-known that a group $G$ is periodic if and only if every subgroup of $G$ contains a minimal subgroup, and the Frattini subgroup $\Phi(G)$ of a group $G$ is equal to the intersection of all maximal subgroups of $G$ if $G$ has at least one maximal subgroup, otherwise it is equal to $G$. Since a minimal (maximal, respectively) subgroup of $G$ is precisely a maximal (primitive, respectively) idempotent of ${\cal K}_1(G)$, it follows that $\Omega_G$ is precisely the set of all maximal idempotents of ${\cal K}_1(G)$. Thus we have immediately

\begin{Lemma}\label{minmax}
\rm Let $G$ be a group, $H\leq G$ and $a\in G$. Then
\begin{itemize}
\item[$(i)$] $Ha\in {\cal K}_1(G)$ is a maximal element if and only if $H\in \Omega_G$;
\item[$(ii)$] $G$ is periodic if and only if there exists $L \in \Omega_G$ such that $K\leq_n L$ for any given $K\in E({\cal K}_1(G))$;
\item[$(iii)$] $H\leq\Phi(G)$ if and only if $K\leq_n H$ for every primitive idempotent $K$ of ${\cal K}_1(G)$.
\end{itemize}
\end{Lemma}

\begin{Lemma}\label{filter}
\rm In the partial order set $({\cal K}_1(G), \leq_n)$, the order filter  $(Ha)^\upharpoonright$ generated by an element $Ha$ is equal to
$ \{Kb \mid 1< K\leq H, \, b\in G, \, ab^{-1} \in H\}$.
In particular, $H^\upharpoonright = {\cal K}_1(H)$.
\end{Lemma}

\begin{Proof}
For any $Kb\in{\cal K}_1(G)$, we have
$$\begin{array}{lll}
Ha \leq_n Kb &\Longleftrightarrow& (\, \exists \, M \in E({{\cal K}_1(G)})) \, Ha= M* Kb \\
&\Longleftrightarrow& (\, \exists \, 1< M \leq G) \, Ha= \langle M, K\rangle b \\
&\Longleftrightarrow&(\, \exists \, 1< M \leq G) \, H= \langle M, K\rangle , \, ab^{-1}\in H \\
&\Longleftrightarrow& K\leq H , \, ab^{-1}\in H. \\
\end{array}$$
\vskip -6pt
This shows that $(Ha)^\upharpoonright = \{Kb \mid 1< K\leq H, \, b\in G, \, ab^{-1} \in H\}$.
\end{Proof}

\begin{Lemma}\label{order}
\rm  Let $G$ be a finite group and $\pi_k$ denote the set of all prime divisors of $\prod_{H\in \Omega_G}|R_{\scriptscriptstyle H}|$, where $R_{\scriptscriptstyle H}$ denotes the ${\cal R}$-class of ${\cal K}_1(G)$ containing $H$.
If $|\pi_k|>1$ (i.e., $G$ is not a $p$-group), then $|G|={\rm l.c.m.}\{|R_{\scriptscriptstyle H}|: H\in \Omega_G\}$. Otherwise, $|\pi_k|=1$, say, $\pi_k=\{p\}$ for some prime $p$ (i.e., $G$ is a $p$-group), then $|G|=p\cdot {\rm l.c.m.}\{|R_{\scriptscriptstyle H}|: H\in \Omega_G\}$.
\end{Lemma}

\begin{Proof}
\rm Suppose that $|G|=p_1^{\alpha_1}p_2^{\alpha_2}\cdots p_n^{\alpha_n}$, where $p_1,\ldots,p_n$ are distinct primes.
Then we know by the Cauchy Theorem (see \cite[1.6.17]{robinson}) that there exist $H_i\in \Omega_G$ with $|H_i|=p_i$ for each $p_i$, $i=1,2\ldots,n$. Hence, we can see from Lemma \ref{kgjiben} that
$$|R_{\scriptscriptstyle H_i}|=|G:H_i|=p_1^{\alpha_1}\cdots p_{i-1}^{\alpha_{i-1}}p_i^{\alpha_i-1}
p_{i+1}^{\alpha_{i+1}}\cdots p_n^{\alpha_n}.$$
This implies that
$$|\pi_k|>1 \iff n>1 \iff \text{ $G$ is not a $p$-group,}$$
and $|G|={\rm l.c.m.}\{|R_{\scriptscriptstyle H}|: H\in \Omega_G\}$ if $|\pi_k|>1$, and $|G|=p\cdot {\rm l.c.m.}\{|R_{\scriptscriptstyle H}|: H\in \Omega_G\}$ if $\pi_k=\{p\}$ for some prime $p$, as required.
\end{Proof}

Lemma \ref{filter} tells us that the filter $H^\upharpoonright$ of ${\cal K}_1(G)$  is an inverse subsemigroup of ${\cal K}_1(G)$ for each $H\in E({\cal K}_1(G))$.
Combining Lemma \ref{order} and \cite[Theorem 1.2.7]{schmi}, we have
\begin{Lemma}\label{cyclic}
\rm  Let $H$ be a finite subgroup of a group $G$. Then $H$ is a cyclic subgroup of order $p_1^{n_1}p_2^{n_2}\cdots p_s^{n_s}$ if and only if $E(H^\upharpoonright)^1$ is a direct product of chains of length $n_1,n_2,\ldots,n_s$, where $p_1,p_2,\ldots,p_s$ are distinct primes.
\end{Lemma}

\begin{Lemma}\label{perm}
\rm  Let $G$ be a group and $H$, $K\in \Omega_G$ with $H\neq K$. Then $H$ and $K$ permute in $G$ if and only if $|R_{ \scriptscriptstyle H}^{\scriptscriptstyle (H*K)^\upharpoonright}|$ is a prime, where $R_{ \scriptscriptstyle H}^{\scriptscriptstyle (H*K)^\upharpoonright}$ denotes the ${\cal R}$-class of the inverse subsemigroup $(H*K)^\upharpoonright$ containing $H$.
\end{Lemma}

\begin{Proof}
We know that $H$ and $K$ permute if and only if $HK=KH$, i.e., $HK=\langle H,K\rangle$ in $G$. The latter one holds if and only if
$|\langle H,K\rangle:H|=|HK:H|$. Noticing that
$$|HK:H|=|K:K\cap H|=|K|/|H\cap K|=|K|.$$
Also, we know that $|K|$ is a prime, since a maximal idempotent of ${\cal K}_1(G)$ is a minimal subgroup of $G$. Thus we can see that $H$ and $K$ permute in $G$ if and only if $|R_{ \scriptscriptstyle H}^{\scriptscriptstyle (H*K)^\upharpoonright}|$ is a prime, since $|R_{ \scriptscriptstyle H}^{\scriptscriptstyle (H*K)^\upharpoonright}|=|\langle H,K\rangle:H|$.
\end{Proof}

To describe the completion of ${\cal K}_1(G)$, we need to study the compatible relations in ${\cal K}(G)$, since the compatible relation on ${\cal K}_1(G)$ is exactly the restriction of the compatible relation on ${\cal K}(G)$ to ${\cal K}_1(G)$.

\begin{Lemma}\label{hakbc}
\rm Let $G$ be a group and $Ha, Kb\in {\cal K}(G)$. Then
$$Ha\thicksim Kb\iff ab^{-1}\in \langle H,K^{ba^{-1}}\rangle\cap \langle H, K\rangle.$$
\end{Lemma}

\begin{Proof}
We can see from Lemma \ref{kgjiben} that
$$\begin{aligned}
Ha\thicksim Kb&\iff Ha*(Kb)^{-1},\,(Ha)^{-1}*Kb\in E({\cal K}_1(G))\\
&\iff \langle H,K^{ba^{-1}}\rangle ab^{-1},\,\langle H^a, K^a\rangle a^{-1}b\leq G\\
&\iff ab^{-1}\in \langle H,K^{ba^{-1}}\rangle,\,a^{-1}b\in  \langle H^a, K^a\rangle\\
&\iff ab^{-1}\in \langle H,K^{ba^{-1}}\rangle\cap \langle H, K\rangle.
\end{aligned}$$
\end{Proof}

\begin{Corollary}\label{hak}
\rm Let $G$ be a group, $a\in G$ and $H, K\leq G$ with $HK=KH$. Then the following statements are equivalent:
\begin{itemize}
\item[$(i)$] $Ha\thicksim K;$
\item[$(ii)$] $a\in HK;$
\item[$(iii)$]$(\exists\, b\in K)\,\, Ha=Hb.$
\end{itemize}
\end{Corollary}

\begin{Proof}
By Lemma \ref{hakbc},  $(i)\Rightarrow (ii)$ and $(ii)\Rightarrow (iii)$ are evident. To show $(iii)\Rightarrow (i)$, notice that
$$Ha*K^{-1}=Hb*K= \langle H,K^{b^{-1}}\rangle b= \langle H,K\rangle b=\langle H,K\rangle \in E({\cal K}_1(G)),$$
$$(Ha)^{-1}*K=H^{b^{-1}}b*K=\langle H^{b^{-1}}, K^{b^{-1}}\rangle b= \langle H^{b^{-1}}, K\rangle b= \langle H^{b^{-1}}, K\rangle \in E({\cal K}_1(G)).$$
Hence $Ha\thicksim K$, as required.
\end{Proof}

As a consequence of  \cite[Lemma 1.3]{lawson1}, we have immediately
\begin{Lemma}\label{rjiaowei1}
\rm Let $S$ be an inverse semigroup. Then $\thicksim \cap \,{\cal L}=\thicksim \cap \,{\cal R}=\Delta_S$.
\end{Lemma}

\begin{Lemma} \label{jiaowei1}
\rm Let $G$ be a group, $\{K_\lambda \}_{\lambda \in \Lambda} = \{K \mid 1 \neq K \leq G\}$ and $A = \{K_\lambda a_\lambda\}_{\lambda \in \Lambda_{\scriptscriptstyle A}}\in C({\cal K}_1(G))$. Then $|A\cap L_{\scriptscriptstyle {K_\lambda}}|\leq 1$ and $|A\cap R_{\scriptscriptstyle {K_\lambda}}|\leq1$ for any $\lambda\in\Lambda$. Moreover,
$$A\in \Sigma({\cal K}_1(G)) \iff
(\forall \lambda \in \Lambda)\,\,|A\cap L_{\scriptscriptstyle {K_\lambda}}|=|A\cap R_{\scriptscriptstyle {K_\lambda}}|=1.$$
\end{Lemma}

\begin{Proof}
It is easy to see from Lemma $\ref{rjiaowei1}$ that $|A\cap L_{\scriptscriptstyle {K_\lambda}}|\leq 1$ and $|A\cap R_{\scriptscriptstyle {K_\lambda}}|\leq1$ for any $\lambda\in\Lambda$.
Also,  Lemma \ref{kgjiben} tells us that $ E({\cal K}_1(G))=\{K_\lambda\}_{\lambda \in \Lambda}$,
and we have from \cite[Chapter 1, Lemma 22]{lawson}
\begin{eqnarray}\label{21}
AA^{-1} & = & \{(K_\lambda a_\lambda)*(K_\lambda a_\lambda)^{-1}\}_{\lambda \in \Lambda_{\scriptscriptstyle A}} = \{K_\lambda\}_{\lambda \in \Lambda_{\scriptscriptstyle A}},
\end{eqnarray}
\begin{eqnarray}\label{22}
A^{-1}A & = &\{(K_\lambda a_\lambda)^{-1}*(K_\lambda a_\lambda)\}_{\lambda \in \Lambda_{\scriptscriptstyle A}} = \{K_\lambda ^{a_\lambda}\}_{\lambda \in \Lambda_{\scriptscriptstyle A}}.
\end{eqnarray}

Suppose that $A\in \Sigma({\cal K}_1(G))$ and $AB = BA = E({\cal K}_1(G))$. Then $B$ is exactly the inverse of $A$ in the inverse semigroup  $C({\cal K}_1(G))$. That is to say, $B = A^{-1}$ and so
$AA^{-1}=A^{-1}A=E({\cal K}_1(G))$. Thus it follows from (\ref{21}) and (\ref{22}), respectively that  $\Lambda_A=\Lambda$ and there exists $K_\nu a_\nu\in A$ such that $K_\nu^{a_\nu} = K_\lambda$ for any $\lambda\in \Lambda$. This implies that $K_\nu a_\nu \, {\cal L} \, K_\lambda$ and $K_\lambda a_\lambda \, {\cal R} \, K_\lambda$ for any $\lambda\in \Lambda$
and so $1 \leq |A\cap L_{\scriptscriptstyle {K_\lambda}}|$, $1 \leq |A\cap R_{\scriptscriptstyle {K_\lambda}}|$ for any $\lambda \in \Lambda$.
Now we have shown that
$$(\forall \lambda \in \Lambda)\,\,|A\cap L_{\scriptscriptstyle {K_\lambda}}|=|A\cap R_{\scriptscriptstyle {K_\lambda}}|=1.$$

Conversely, suppose that $|A\cap L_{\scriptscriptstyle {K_\lambda}}|=|A\cap R_{\scriptscriptstyle {K_\lambda}}|=1$ for any $\lambda \in \Lambda$. Then there exist $K_\mu a_\mu, \, K_\nu a_\nu\in A$ such that $K_\mu a_\mu \, {\cal L} \, K_\lambda$ and $K_\nu a_\nu \, {\cal R} \, K_\lambda$ for any $\lambda\in \Lambda$. That is to say, $K_\lambda=K_\nu\in AA^{-1}$ and $K_\lambda=K_\mu^{a_\mu} \in A^{-1}A$. This shows that $AA^{-1}=A^{-1}A=E({\cal K}_1(G))$ and so $A\in \Sigma({\cal K}_1(G))$, as required.
\end{Proof}

\section{The unit group $\Sigma({\cal K}_1(G))$}
In this section we shall investigate the unit group $\Sigma({\cal K}_1(G))$ of the completion $C({\cal K}_1(G))$ of ${\cal K}_1(G)$.
Recall that $\{K_\lambda \}_{\lambda \in \Lambda} = \{K \mid 1 \neq K \leq G\}$ and $\Omega_G = \{H_i\}_{i\in I}$ denotes the set of all minimal subgroups of $G$ and
$${\cal A}_G=\{\{H_ia_i\}_{i\in I}\mid (\forall \,j, k\in I)\,H_ja_j\thicksim H_ka_k \text{ and } |\{H_ia_i\}_{i\in I}\cap L_{\scriptscriptstyle {H_j}}|=1\}.$$
It is well known that ${\cal K}_1(G)$ is an inverse semigroup and a partially ordered set under the natural order. For any $A\subseteq {\cal K}_1(G)$ we denote by $[A]$ the order ideal generated by $A$, and by $\overline A$  the set of all maximal elements in $A$.
In particular, for any $U=\{K_{\lambda}a_{\lambda}\}_{\lambda\in \Lambda}\in \Sigma({\cal K}_1(G))$, $\overline U=\{H_ia_i\}_{i\in I}$, and $[\,\overline{U}\,]=U$.

Firstly, we have

\begin{Lemma}\label{tonggou}
\rm  ${\cal A}_G$ equipped with the multiplication $A\cdot B=\overline {[A][B]}$
becomes a group, and the mapping $\varphi$ from ${\cal A}_G$ to $ \Sigma({\cal K}_1(G))$ defined by $\varphi (A) = [A]$ is an isomorphism.
\end{Lemma}

\begin{Proof}
\rm Let $A=\{H_ia_i\}_{i\in I}\in {\cal A}_G$. It follows immediately from \cite[Section 1.4, Lemma 14 ]{lawson} that $[A]\in C({\cal K}_1(G))$ and so $|[A]\cap L_{\scriptscriptstyle {K}}|\leq 1$ and $|[A]\cap R_{\scriptscriptstyle {K}}|\leq1$ for any $1\neq K\leq G$ by Lemma \ref{jiaowei1}. Also, there exists $H_s\in\Omega_G$ such that $H_s\leq K$, since $G$ is periodic. It follows that there exists $H_ta_t\in A$ such that $H_ta_t \, {\cal L} \, H_s$, since $|\{H_ia_i\}_{i\in I}\cap L_{\scriptscriptstyle {H_s}}|=1$. That is to say, $H_t^{a_t}=H_s\leq K$. This implies that $Ka_s\leq_n H_sa_s$ and $K^{a_t^{-1}}a_t\leq_n H_ta_t$ and so $Ka_s, \, K^{a_t^{-1}}a_t\in [A]$. Noticing that $Ka_s \, {\cal R} \, K$ and $K^{a_t^{-1}}a_t \, {\cal L} \, K$, we have that $|[A]\cap L_{\scriptscriptstyle {K}}|= 1$ and $|[A]\cap R_{\scriptscriptstyle {K}}|=1$ for any $1\neq K\leq G$. Thus
$[A]\in \Sigma({\cal K}_1(G))$ by Lemma \ref{jiaowei1}. We have shown that $\varphi$ is well-defined. Also $\varphi$ is injective.

Suppose that $U\in \Sigma({\cal K}_1(G))$ and write $A=\overline{U}$. Then $|U\cap L_{\scriptscriptstyle {K}}|=|U\cap R_{\scriptscriptstyle {K}}|=1$ for any $1\neq K\leq G$ by Lemma \ref{jiaowei1}. It is clear that $A$ is compatible and we may write $A=\{H_ia_i\}_{i\in I}$ by Lemma \ref{minmax}. Also, for any $H_j\in\Omega_G$, we have $|A\cap L_{\scriptscriptstyle {H_j}}|\leq1$ by Lemma \ref{rjiaowei1} and there exists $Ha\in U$ such that $Ha \, {\cal L} \, H_j$. That is to say, $H^a=H_j$ and so $H=H_j^{a^{-1}}\in\Omega_G$. It follows from Lemma \ref{minmax} that $Ha\in A$ and so $A\in {\cal A}_G$. Noticing that $\varphi(A)=[A]=[\overline{U}]=U$, we have shown that $\varphi$ is surjective and the multiplication $\cdot$ on ${\cal A}_G$ defined above is well-defined. Also, we have that
$$\varphi(A\cdot B)=[A\cdot B]=\left[\overline{[A] [B]}\right]= [A] [B]=\varphi(A) \varphi(B),$$
for all $A,B\in{\cal A}_G$, since $[A] [B]\in \Sigma({\cal K}_1(G))$.  Hence $\varphi$ is bijective and preserves the multiplications, and so $({\cal A}_G,\cdot)$ is a group, since $\Sigma({\cal K}_1(G))$ is a group, as required.
\end{Proof}

Suppose that $A=\{H_ia_i\}_{i\in I}$, $B=\{H_ib_i\}_{i\in I}\in {\cal A}_G$, and $A\cdot B=\{H_i c_i\}_{i\in I},$
then for any $k\in I$, there exist $i,j\in I$ such that
$$H_kc_k=H_ia_i* H_ja_j=\langle H_i,\,H_j^{a_i^{-1}}\rangle a_ib_j.$$
Hence $\langle H_i,\,H_j^{a_i^{-1}}\rangle = H_k$, and so $H_i=H_j^{a_i^{-1}}=H_k$, thus $H_j=H_i^{a_i}$.
Further,  defining  a transformation $\tau_A$ on $I$  by $i^{\tau_A}=j$ whenever $H_j=H_i^{a_i}$, we have
$$A\cdot B=\{H_ia_ib_{i^{\tau_A}}\}_{i\in I}.$$
In fact, we can prove that $\tau_A$ is a permutation on $I$. It is easy to see that $\tau_A$ is well- defined, since $H_ia_i=H_ia^{\prime}_i$ implies that $H_i^{a_i}=H_i^{a'_i}$. Suppose that $i^{\tau_A}=j^{\tau_A}$. Then $H_i^{a_i}=H_j^{a_j}$ and so $H_ia_i \, {\cal L} \, H_ja_j$.
Noticing that $H_ia_i\thicksim H_ja_j$, we have that $H_ia_i=H_ja_j$ by Lemma \ref{rjiaowei1}. Thus $i=j$ and so $\tau_A$ is injective. Also, there exists $H_ia_i\in A$ such that $H_ia_i \, {\cal L} \, H_j$ for any $j\in I$, since $|A\cap L_{\scriptscriptstyle {H_j}}|=1$. It follows that $H_i^{a_i}=H_j$ and so $i^{\tau_A}=j$. This shows that $\tau_A$ is surjective and so it is a permutation on $I$.

Define a map $\tau$ from ${\cal A}_G$ to $Sym(I)$ by $\tau(A)=\tau_A$. Then for each $k\in I$, one has
$$H_{k^{\tau(A)\tau(B)}}=H_{(k^{\tau_A})^{\tau_B}}=(H_{k^{\tau_A}})^{b_{k^{\tau_A}}}=H_k^{a_kb_{k^{\tau_A}}}
=H_{k^{\tau(A\cdot B)}},$$
and so $k^{\tau(A)\tau(B)}=k^{\tau(A\cdot B)}$. This shows that
$$\tau(A)\tau(B)=\tau(A\cdot B).$$ That is to say, $\tau$ is a homomorphism.

\noindent{\bf Remark\,} Let $\{\Omega_t\}_{t\in T}$ be all orbits of $\Omega_G$ under the conjugation action of $G$ on $\Omega_G$ and
$X_t=\{k\in I\mid H_k\in \Omega_t\}$ for any $t\in T$. Then
$$\tau({\cal A}_G)\leq {\textstyle\prod}_{t\in T}\, Sym(X_t),$$
since $H_{i^{\tau_A}}$ is conjugate to $H_{j^{\tau_A}}$ whenever $H_i$ is conjugate to $H_j$.

Summarizing the above results, we have
\begin{Proposition}
\rm For any $A=\{H_ia_i\}_{i\in I}$, $B=\{H_ib_i\}_{i\in I}\in {\cal A}_G$,
$$A\cdot B=\{H_ia_ib_{i^{\tau_A}}\}_{i\in I},$$
 where the transformation $\tau_A$ defined by $i^{\tau_A}=j$ whenever $H_j=H_i^{a_i}$ is a permutation on $I$ and the map $\tau$ from ${\cal A}_G$ to $Sym(I)$ defined by $\tau(A)=\tau_A$ is a group homomorphism.

\end{Proposition}

Next, we shall investigate the unit group $\Sigma({\cal K}_1(G))$ when $G$ is a nontrivial periodic group satisfying $H_iH_j = H_jH_i$ for any $H_i, \,H_j\in \Omega_G$ and $|\Omega_G|\geq 2$. Write $\Omega_p(G)$ to denote the subgroup of $G$ generated by all minimal subgroups of $G$ of order $p$, i.e., $\Omega_p(G)=\langle H_i\in \Omega_G \mid |H_i|=p\rangle=\langle x\in G\mid x^p=1\rangle$. Let $\Omega(G)=\prod_{p\in\mathbb{P}}\Omega_p(G)$, where $\mathbb{P}$ denotes the set of all primes. Then we have

\begin{Lemma}\label{minper}
\rm The following statements are equivalent:
\begin{itemize}
\item[$(i)$] $H_iH_j = H_jH_i$ for any $H_i, \,H_j\in \Omega_G$;
\item[$(ii)$] $\Omega_p(G)$ is either a trivial subgroup or an elementary abelian $p$-subgroup of $G$ for each prime $p$;
\item[$(iii)$] $\Omega(G)$ is an abelian subgroup of $G$.
\end{itemize}
\end{Lemma}

\begin{Proof}
Firstly, we prove that $(i)\Rightarrow(ii)$.
Suppose that $p$ is a prime and $\Omega_p(G)$ is not trivial. It is sufficient to show that $\Omega_p(G)$ is abelian. Suppose $a, b\in G$ with $|a|=|b|= p$. Then $\langle a\rangle=H_i$, $\langle b\rangle=H_j$ for some $H_i$, $H_j\in \Omega_G$, and $H_iH_j\leq G$. Also, we can see that $|H_iH_j|=p$ or $|H_iH_j|=p^2$,
since
$$|H_iH_j|=|H_i||H_j|/|H_i\cap H_j|=p^2/|H_i\cap H_j|.$$
This implies that  $H_iH_j$ is an abelian group by \cite[1.6.15]{robinson} and so $ab = ba$, as required.

Now, we prove that $(ii)\Rightarrow(iii)$. It is easy to see that $\Omega_p(G)\unlhd G$ for any prime $p$. This implies that $\Omega(G)$ is the direct product of $\Omega_p(G)$, $p\in\mathbb{P}$. Thus $\Omega(G)$ is abelian.

It is easy to see $(iii)\Rightarrow(i)$, since $\Omega(G)$ contains all minimal subgroups of $G$. This completes the proof.
\end{Proof}

\begin{Lemma}\label{ag}
\rm  Let $Ag=\{H_i a_ig\}_{i\in I}$ for any $A=\{H_i a_i\}_{i\in I}\in {\cal A}_G$ and $g\in G$. Then $Ag \in {\cal A}_G$.
\end{Lemma}

\begin{Proof}
It is easy to see from Lemma \ref{hakbc} that $H_ia_ig\thicksim H_ja_jg$ for any $i, \, j \in I$, since $H_ia_i\thicksim H_ja_j$. Also, there exists $H_ia_i \in A$ such that $H_ia_i \, {\cal L} \, H_j^{g^{-1}}$ for any $j\in I$, since $|A\cap L_{\scriptscriptstyle{H_j^{g^{-1}}}}|=1$. That is to say, $H_i^{a_i}=H_j^{g^{-1}}$, i.e., $H_i^{a_ig}=H_j$. This implies that $H_i{a_ig}\, {\cal L} \,H_j$ and so $|Ag \cap L_{\scriptscriptstyle {H_j}}|=1$ for any $j\in I$. Thus $Ag\in {\cal A}_G$, as required.
\end{Proof}

In the remainder of this section, we always assume that $G$ is a nontrivial periodic group satisfying $H_iH_j = H_jH_i$ for any $H_i, \,H_j\in \Omega_G$ and $|\Omega_G|\geq 2$.
Given $H_r \in \Omega_G$.  Taken a right transversal ${\cal{C}}$ of $H_r$ in $G$. Write ${\cal A}(H_r)$ to denote the subset of ${\cal A}_G$ consisting of all $\{H_i a_i\}_{i\in I}\in{\cal A}_G$
such that $ a_r \in H_r $.
Define a map
$\psi$ from  ${\cal A}(H_r)\times {\cal C}$ to ${\cal A}_G$ by $\psi(A, c_{\scriptscriptstyle\ell})= Ac_{\scriptscriptstyle\ell}$.
Then we have

\begin{Lemma}\label{t-1}
\rm The mapping $\psi$ defined above is a bijection. In particular,
$$|{\cal A}_G|=|{\cal A}(H_r)||G:H_r|.$$
\end{Lemma}

\begin{Proof}
We can see from Lemma \ref{ag} that $\psi$ is well-defined.
Suppose that $Ac_{\scriptscriptstyle\ell}=Bc_k$ for some $A,\, B\in {\cal A}(H_r)$ and $c_{\scriptscriptstyle\ell}, \, c_k\in {\cal C}$. Then $H_rc_{\scriptscriptstyle\ell}=H_rc_k$ and so $c_{\scriptscriptstyle\ell}=c_k$. This implies that $$B=(Bc_k)c_k^{-1}=(Ac_{\scriptscriptstyle\ell})c_{\scriptscriptstyle\ell}^{-1}=A$$
and so $\psi$ is injective.

For each $A=\{H_ia_i\}_{i\in I}\in{\cal A}_G$, there exists $c_{\scriptscriptstyle\ell}\in\mathcal C$ such that $H_r a_r = H_r c_{\scriptscriptstyle\ell}$. Let $A'=Ac_{\scriptscriptstyle\ell}^{-1}$.
Then $A'\in {\cal A}(H_r)$ by Lemma \ref{ag} and $\varphi(A',c_{\scriptscriptstyle\ell})=A$. This shows that $\psi$ is bijective, and so
$$|{\cal A}_G|=|{\cal A}(H_r)\times {\cal C}|=|{\cal A}(H_r)| |{\cal C}|=|{\cal A}(H_r)||G:H_r|,$$
as required.
\end{Proof}

Define a simple graph $\Gamma(H_r)$ corresponding to $H_r$ by $\Gamma(H_r)=(V_{\scriptscriptstyle{H_r}},E_{\scriptscriptstyle{H_r}})$, where $V_{\scriptscriptstyle{H_r}}=\Omega_G\setminus \{H_r\}$ and
$$E_{\scriptscriptstyle{H_r}}=\{(H_i,H_j)\mid H_i\neq H_j,\,\,H_iH_j\cap H_r=1\}.$$
Write $\{\Gamma_j\}_{j\in J}$ to denote the set of all connected components of $\Gamma(H_r)$.
We have

\begin{Lemma}\label{graph}
\rm The followings are true:
\begin{itemize}
\item[$(i)$] ${\cal A}(H_r)=\left\{\bigcup_{j\in J}\Gamma_ja_j\bigcup \{H_r\}\mid (\forall j\in J)\, a_j\in H_r\right\},$
\item[$(ii)$] $|{\cal A}(H_r)|=|H_r|^{|J|}.$
\end{itemize}
\end{Lemma}

\begin{Proof}
\rm We can see from Corollary \ref{hak} that
$$\begin{aligned}
H_ia_i\thicksim H_r\iff&(\exists\, b_i\in H_r)\,\, H_ia_i=H_ib_i,\\
H_ia_i\thicksim H_ja_j\iff& a_ia_j^{-1}\in H_iH_j\iff b_ib_j^{-1}\in H_iH_j,
\end{aligned}$$
for any given $A=\{H_ia_i\}_{i\in I}\in {\cal A}(H_r)$.
It follows that
$$A=\{H_ib_i\}_{i\in I\setminus\{t\}}\cup\{H_r\},$$
where $b_i\in H_r$, $b_ib_j^{-1}\in H_iH_j$ for all $i,j\in I\setminus\{r\}$.

On the other hand, $b_i=b_j$ if $(H_i, H_j)\in E_{\scriptscriptstyle{H_r}}$, since $b_ib_j^{-1}\in H_iH_j\cap H_r=1$. This implies that $b_i=b_j$ if $H_i, H_j\in \Gamma_k$ for some $k\in J$.  This shows that
$$A=\cup\{\Gamma_jc_j\}_{j\in J}\cup \{H_r\}$$
for some $c_j\in H_r,\,j\in J$.

Conversely, suppose that $A=\cup\{\Gamma_ja_j\}_{j\in J}\cup \{H_r\}$, where $a_j\in H_r$ for any $j\in J$.
Then, it is easy to see that $A$ is a compatible subset of ${\cal K}_1(G)$. Also, we can see from Lemma \ref{minper} that $H_i^{a_j}=H_i$, i.e., $H_ia_j  \, {\cal L} \, H_i$ for any $H_ia_j\in A$. Hence $|A \cap L_{\scriptscriptstyle {H_i}}|=1$ for any $H_i\in\Omega_G$ and so $A\in  {\cal A}(H_r)$. This shows that
$${\cal A}(H_r)=\{\cup\{\Gamma_ja_j\}_{j\in J}\cup \{H_r\}\mid (\forall j\in J)\, a_j\in H_r\},$$
and so $|{\cal A}(H_r)|=|H_r|^{|J|}$, as required.
\end{Proof}

By Lemma \ref{tonggou}, Lemma \ref{t-1} and Lemma \ref{graph} we immediately have

\begin{Theorem}\label{|ag|}
$|\Sigma({\cal K}_1(G))|=|{\cal A}_G|=|G||H_r|^{|J|-1}.$
\end{Theorem}

\rm Next, we shall show that $\Gamma(H_r)$ is connected if $G$ is not a $p$-group for any  prime $p$ or $G$ is a $p$-group with $|\Omega_G|> p^2$ for some prime $p$, otherwise it is an empty graph.

\begin{Theorem}\label{pq}
\rm If $G$ is not a $p$-group for any prime $p$ and $H\in V_{\scriptscriptstyle H_r}$ with $|H|\neq |H_r|$, then $(H_i,H)\in E_{\scriptscriptstyle H_r}$ for all $H_i\in V_{\scriptscriptstyle H_r}$ with $H_i\neq H$, and $|J|=1$.
\end{Theorem}

\begin{Proof}
\rm Suppose that $|H|=p$ and $|H_r|=q$.  If $(H_i,H)\notin E_{\scriptscriptstyle H_r}$ for some $H_i\in E_{\scriptscriptstyle H_r}$ with $H_i\neq H$, then $H_iH\cap H_r=H_r$, and so $H_r\leq H_iH$. Thus we have
$$p |H_i| = |H||H_i|/|H \cap H_i| = |H H_i| = |H H_i : H_r||H_r| = q |H H_i : H_r|.$$
This implies that $|H_i|=q$, since $p,\,q $ and $|H_i|$ are prime numbers and $p \neq q$.
We can now see that $q^2 = |H_i||H_r| = |H_iH_r|$ and $ pq =|H||H_i| = |HH_i|$.
Again by $H_r\leq HH_i$ and $H_i\leq HH_i$, we have that $H_iH_r\leq HH_i$ and so $q^2$ divides $pq$,
a contradiction. This shows that $(H,H_i)\in E_{\scriptscriptstyle H_r}$ for all $H_i\in V_{\scriptscriptstyle H_r}$ with $H_i\neq H$, and so $|J|=1$, as required.
\end{Proof}

Now, we shall consider the case where $G$ is a $p$-group. In this case $\Omega(G)=\Omega_p(G)$ is an elementary abelian $p$-subgroup of $G$ (see Lemma \ref{minper}).
It is well known that the elementary abelian group $\Omega(G)$ equipped with the scalar multiplication defined by $\overline k\circ \alpha=k\alpha$  becomes a vector space over field $\mathbb Z_p$. Also, $H_i$ is a $1$-dimensional subspace of $\Omega(G)$ for every $H_i\in \Omega_G$, since $H_i=\{\overline{k}a_i\mid \overline{k}\in\mathbb Z_p\}$ for some $a_i\in \Omega(G)$.

To study $\Gamma(H_r)$ for the $p$-groups, the following two lemmas are needed for us.
\begin{Lemma}\label{zikongjiantu}
\rm Let $W$ be a vector space over field $\mathbb F$ such that $\dim W\geqslant 2$ and $v\in W\setminus \{0\}$. Define a graph $\Gamma_v=(V_{\scriptscriptstyle \Gamma_v},E_{\scriptscriptstyle \Gamma_v})$, where
$V_{\scriptscriptstyle \Gamma_v}=\{L(\alpha)\mid \alpha\in W\setminus\{0\}\}\setminus\{L(v)\}$ and
$$
E_{\scriptscriptstyle \Gamma_v} =
\{(L(\alpha),L(\beta))\mid L(\alpha)\neq L(\beta),\,\, v\notin L(\alpha,\beta)\}.
$$
Then
$$\ell_{\scriptscriptstyle\Gamma_v}=\begin{cases}
|V_{\scriptscriptstyle \Gamma_v}|, &\dim W=2,\\
1,  &\dim W\geqslant 3,
\end{cases}$$
where $\ell_{\scriptscriptstyle\Gamma_v}$ denotes the cardinal number of the set of all connected components of $\Gamma_v$.
\end{Lemma}

\begin{Proof}
Let $L(\alpha),\, L(\beta)\in V_{\scriptscriptstyle \Gamma_v}$ with $L(\alpha)\neq L(\beta)$. Then $\dim L(\alpha,\beta)=2$. If $\dim W=2$, then $W=L(\alpha,\beta)$, and so $v\in L(\alpha,\beta)$. This implies that $\left(L(\alpha),L(\beta)\right)\notin E_{\scriptscriptstyle \Gamma_v}$ and $\ell_{\scriptscriptstyle\Gamma_v}= |V_{\scriptscriptstyle \Gamma_v}|$. Otherwise, $\dim W\geqslant 3$. Suppose that $\left(L(\alpha),L(\beta)\right)\notin E_{\scriptscriptstyle \Gamma_v}$. Then $v\in L(\alpha, \beta)$. Thus $L(\alpha, \beta,v)=L(\alpha,\beta)\neq W$, since $\dim W\geqslant 3$.
We can show that
$$(L(\alpha),L(\theta)),\,(L(\theta),L(\beta))\in E_{\scriptscriptstyle \Gamma_v}$$
for any $\theta\in W\setminus L(\alpha, \beta,v)$.
In fact,  if $(L(\alpha),L(\theta))\notin E_{\scriptscriptstyle \Gamma_v}$, then $v\in L(\alpha, \theta)$, and so $v=k_1\alpha+k_2\theta$ for some $k_1,k_2\in \mathbb F$ with $k_2\neq 0$, since $L(v)\neq L(\alpha)$. This implies that
$$\theta=1/k_2(v-k_1\alpha)\in L(\alpha,\beta,v),$$
a contradiction. Thus $(L(\alpha),L(\theta))\in E_{\scriptscriptstyle \Gamma_v}$.
Dually, we can prove that $(L(\theta),L(\beta))\in E_{\scriptscriptstyle \Gamma_v}$. This shows that $\ell_{\scriptscriptstyle\Gamma_v}=1$, as required.
\end{Proof}

The following result is well-known.

\begin{Lemma}\label{yiwiezikongjian}
\rm Let $W$ be a $m$-dimensional vector space over field $\mathbb Z_p$, where $p$ is a prime number.
Then $W$ has $(p^m-1)/(p-1)$ $1$-dimensional subspaces, i.e.,
$$ \left|\{L(\alpha)\mid \alpha\in W\setminus\{0\}\}\right| = (p^m-1)/(p-1).$$
\end{Lemma}

By Lemma $\ref{minper}$, Lemma $\ref{zikongjiantu}$ and Lemma $\ref{yiwiezikongjian}$ we immediately have

\begin{Theorem}\label{pm}
\rm If $G$ is a $p$-group, then $\Gamma(H_r)$ is empty when $|\Omega(G)|=p^2$, and $\Gamma(H_r)$ is connected when $|\Omega(G)|>p^2$.
\end{Theorem}

As a consequence, by Theorem \ref{|ag|} we have
\begin{Corollary}
\rm $|\Sigma({\cal K}_1(G))|=|G|$ if $G$ is not a $p$-group for any  prime $p$ or $G$ is a $p$-group with $|\Omega_G|> p^2$ for some prime $p$.
\end{Corollary}
\section{Groups with isomorphic coset semigroups}

In this section we shall directly confront the question: which group $G$ and its coset semigroup ${\cal K}_1(G)$ are uniquely determined by each other, up to isomorphism.

The following shows that a group $G$ is embedded to the unit group $\Sigma({\cal K}_1(G))$ of $C({\cal K}_1(G))$ if the intersection of all non-trivial subgroups of $G$ is trivial.

\begin{Lemma}\label{embedding}
\rm  Let $G$ be a group and $\{K_\lambda\}_{\lambda\in \Lambda}$ the set of all non-trivial subgroups of $G$. Define a map $\eta$ from $G$ to $\Sigma({\cal K}_1(G))$  by $\eta(a)=\{K_\lambda a\}_{\lambda\in \Lambda}$. Then $\eta$ is a group homomorphism, and $\ker\eta=\bigcap\{K_\lambda\}_{\lambda\in \Lambda}$.
\end{Lemma}

\begin{Proof}
\rm By virtue of Lemma \ref{hakbc} and Lemma \ref{jiaowei1}, $\eta(a)\in \Sigma({\cal K}_1(G))$ for each $a\in G$.  That is to say, $\eta$ is well-defined. For any $a,b\in G$, we have
$$
\eta(a)\eta(b)=\{K_\lambda a\}_{\lambda\in \Lambda}\{K_\lambda b\}_{\lambda\in \Lambda}
 =\{K_\lambda a* K_\mu b\}_{\lambda,\,\mu\in\Lambda}
 =\{\langle K_\lambda,K_\mu^{a^{-1}}\rangle ab\}_{\lambda,\,\mu\in\Lambda}
\subseteq \eta(ab).$$
On the other hand, $K_\lambda ab=K_\lambda a*K_\lambda^a b\in\eta(a)\eta(b)$ for any $K_\lambda ab\in \eta(ab)$. This implies that $ \eta(a)\eta(b) \supseteq\eta(ab)$. Thus we have shown that $\eta(a)\eta(b)=\eta(ab)$, and so $\eta$ is a group homomorphism.

Also, for any $a\in G$, we have
$$a\in \ker\eta\iff \eta(a)=\{K_\lambda\}_{\lambda\in \Lambda}\iff (\forall \,\lambda\in\Lambda)\,\,a\in K_\lambda\iff a\in \textstyle\bigcap\{K_\lambda\}_{\lambda\in \Lambda},$$
since the identity of $\Sigma({\cal K}_1(G))$ is $\{K_\lambda\}_{\lambda\in \Lambda}$.
This shows that $\ker\eta=\bigcap\{K_\lambda\}_{\lambda\in \Lambda}$, as required.
\end{Proof}

\begin{Lemma}\label{surj}
\rm  Let $G$ be a periodic group such that $H_iH_j=H_jH_i$ for all $H_i, H_j\in \Omega_G$ and $|\Omega_G|\geq2$. If $\eta$ is as defined in Lemma \ref{embedding}, then $\eta$ is an isomorphism from $G$ to $\Sigma({\cal K}_1(G))$ if and only if  $|J|=1$.
\end{Lemma}

\begin{Proof}
\rm  Suppose that $|J|=1$. Then we can see from Lemma \ref{tonggou}, Lemma \ref{t-1} and Lemma \ref{graph} that $\eta$ is surjective. Also, Theorem \ref{pq} and Theorem \ref{pm} tell us that $|J|=1$ if and only if $G$ is not a $p$-group for any prime $p$ or a $p$-group with $|\Omega(G)|>p^2$ for some prime $p$. Hence it follows from Lemma \ref{embedding} that $\ker\eta=1$.  This shows that  $\eta$ is injective, and so $\eta$ is an isomorphism.

Conversely, suppose that $\eta$ is an isomorphism. If $|J|>1$, then by Lemma \ref{t-1} we may take
$$A=\cup\{\Gamma_ja_0\}_{j\in J,\,j\neq i}\cup\Gamma_ia_i\cup \{H_r\}\in {\cal A}(H_r),$$
where $a_i,a_0\in H_r$ and $a_i\neq a_0$.
Also, it follows from  Lemma \ref{tonggou} that there exists $a\in G$ such that $[A]=\eta(a)$, since $\eta$ is surjective.
In particular, we have that $H_ia_i=H_ia$, $H_ja_0=H_ja$  and $H_ra=H_r$ for any $H_i\in\Gamma_i$ and $H_j\in\Gamma_j$, where $j\neq i$.
This implies that
$$a\in H_r\cap H_ia_i=H_ra_i\cap H_ia_i=(H_r\cap H_i)a_i=\{a_i\},$$
and so $a=a_i$. Dually, we can show that $a=a_0$. That is to say, $a_i=a_0$, a contradiction.
This shows that $|J|=1$, as required.
\end{Proof}

The following tells us that for the periodic groups whose minimal subgroups permute with each other, except some special $p$-groups, $G$ and ${\cal K}_1(G)$ are uniquely determined by each other, up to isomorphism.

\begin{Theorem}\label{isom1}
\rm  Let $G$ be a periodic group such that $H_iH_j=H_jH_i$ for all $H_i, H_j\in \Omega_G$ and $|\Omega_G|\geq2$. If $G$ is not a $p$-group for any prime $p$ or $G$ is a $p$-group with $|\Omega(G)|>p^2$ for some prime $p$, then for any group $\widetilde{G}$,
${\cal K}_1(G)\cong {\cal K}_1(\widetilde{G})$ if and only if $G\cong \widetilde{G}$.
\end{Theorem}

\begin{Proof}
\rm    It is evident that $G\cong \widetilde{G}$ implies that  ${\cal K}_1(G)\cong {\cal K}_1(\widetilde{G})$. Conversely, suppose that ${\cal K}_1(G)\cong {\cal K}_1(\widetilde{G})$. Then we have that $\Sigma({\cal K}_1(G))\cong G$ by Theorem \ref{pq}, Theorem \ref{pm} and Lemma \ref{surj}.
Also, we can see from Lemma \ref{minmax} and Lemma \ref{perm} that $\widetilde{G}$ is a periodic group whose minimal subgroups permute with each other. Moreover, we have that $\Gamma(H_r)\cong \Gamma(\widetilde{H_r})$, where $\widetilde{H_r}$ is the image of $H_r$ under the isomorphism from ${\cal K}_1(G)$ to $ {\cal K}_1(\widetilde{G})$, since $(H_i,H_j)\in E_{\scriptscriptstyle{H_r}}$ if and only if $(H_i* H_j)\vee H_r$ is not exist in ${\cal K}_1(G)$. In particular,  $\Gamma(\widetilde{H_r})$ is connected and so $\Sigma({\cal K}_1(\widetilde{G}))\cong \widetilde{G}$ by Lemma \ref{surj}.

Now, ${\cal K}_1(G)\cong {\cal K}_1(\widetilde{G})$ implies that $\Sigma({\cal K}_1(G))\cong \Sigma({\cal K}_1(\widetilde{G}))$ by \cite[Corollary of 1.17]{schei2}. We then have that $G\cong \widetilde{G}$, as required.
\end{Proof}

To answer the above question for the $p$-groups $G$ with $|\Omega(G)| \leq  p^2$ (i.e., $|\Omega(G)|=p$ or $|\Omega(G)|=p^2$), we need the followings.

\begin{Lemma}\label{abelian}
\rm  Let $G$ be a group. Then ${\cal K}_1(G)$ is abelian if and only if $G$ is either an abelian group or else the quaternion group $Q_8$.
\end{Lemma}

\begin{Proof}
\rm It is easy to see that ${\cal K}_1(G)$ is abelian if $G$ is an abelian group. Also, ${\cal K}_1(Q_8)$ is abelian.
In fact, noticing that $Q_8$ is a Dedekind group and $Ha\in Z(Q_8/H)$ for any $Ha\in {\cal K}_1(Q_8)$, we have
$$Ha*Kb=HKab=HKa*HKb=HKb*HKa=HKba=KHba=Kb*Ha,$$
for any $Ha,\,Kb\in {\cal K}_1(Q_8)$.

Conversely, suppose that ${\cal K}_1(G)$ is abelian for some non-abelian group $G$. Then every idempotent of ${\cal K}_1(G)$ is central. Thus, by Lemma \ref{kgjiben}, every subgroup of $G$ is normal, that is, $G$ is a Dedekind group. Then we have $G\cong Q_8\times A$ by \cite[Theorem 5.3.8]{robinson}, where $A$ is an abelian group. If $A\neq 1$, then by Lemma \ref{kgjiben}, $H_{\scriptscriptstyle A}$, the ${\cal K}$-class of ${\cal K}_1(G)$ containing $H$, is equal to $G/A\cong Q_8$, which is a non-abelian subgroup of ${\cal K}_1(G)$, a contradiction. This implies that $A=1$ and so $G\cong Q_8$, as required.
\end{Proof}

\begin{Theorem}\label{omega1}
\rm  Let $G$ be a $p$-group with $|\Omega(G)|=p$. Write $\Omega(G)=K$. If ${\cal K}_1(G)\cong {\cal K}_1(\widetilde{G})$ for some group $\widetilde{G}$, then $G/K\cong \widetilde{G}/\widetilde{K}$, where $\widetilde{K}$ is the unique minimal subgroup of $\widetilde{G}$. Further, $G\cong \widetilde{G}$ if $G$ is a finite group of composite order.
\end{Theorem}

\begin{Proof}
\rm It is easy to see that $\Omega_G=\{K\}$. Suppose that $\varphi$ is an isomorphism from ${\cal K}_1(G)$ to ${\cal K}_1(\widetilde{G})$. Then we have $\Omega_{\widetilde{G}}=\{ \widetilde{K}\}$, where $\widetilde{K}=\varphi(K)$, since $\varphi$ preserves the natural partial order.
Also, it is easy to see that $K$ and $\widetilde{K}$ are respectively the identities of ${\cal K}_1(G)$ and
${\cal K}_1(\widetilde{G})$.
This implies that their groups of units $U({\cal K}_1(G))=H^{\scriptscriptstyle {\cal K}_1(G)}_{\scriptscriptstyle K}$ and $U({\cal K}_1(\widetilde G))=H^{\scriptscriptstyle {\cal K}_1(\widetilde G)}_{\scriptscriptstyle \widetilde K}$ are isomorphic.
On the other hand, it follows from Lemma \ref{kgjiben} that $H^{\scriptscriptstyle {\cal K}_1(G)}_{\scriptscriptstyle K}=G/K$  and $H^{\scriptscriptstyle {\cal K}_1( \widetilde G)}_{\scriptscriptstyle \widetilde K}=\widetilde{G}/\widetilde{K}$. This shows that $G/K\cong \widetilde{G}/\widetilde{K}$.

Now, assume that $G$ is a finite group of composite order.. Then $G>H$, and by Lemma \ref{order}, we have $|\widetilde{G}|=|G|$.
Since $|\Omega_G|=|\Omega_{\widetilde{G}}|=1$, $G$ is either a cyclic group of prime power order or else a generalized quaternion group by \cite[5.3.6]{robinson} and by the Sylow Theorem. Also, the same conclusion holds for $\widetilde{G}$.

Suppose first that $G$ is a cyclic group of prime power order. Then ${\cal K}_1(G)$ is abelian and so is ${\cal K}_1(\widetilde{G})$. Hence we can see from Lemma \ref{abelian} that $\widetilde{G}\cong G$ if $|G|\neq8$. Suppose that $|G|=8$ and $\widetilde{G}\ncong G$. Then $\widetilde{G}\cong Q_8$ and so ${\cal K}_1(\widetilde{G})$ contains three primitive idempotents, i.e., maximal subgroups of $\widetilde{G}$. However, $G$ contains exactly one primitive idempotent, a contradiction. Hence $\widetilde{G}\cong G$.

Suppose now that $G$ is a generalized quaternion group. If $|G|=8$ then we can see that $\widetilde{G}\cong G$ by previous proof. Hence we assume that $|G|>8$. Then ${\cal K}_1(G)$ is non-abelian and so is ${\cal K}_1(\tilde{G})$. This implies that $\widetilde{G}$ is also a generalized quaternion group. Since $|G|=|\widetilde{G}|$, we get $G\cong \widetilde{G}$, as required.
\end{Proof}

In the remainder of this section, we concentrate on the case of finite groups. We shall show that for finite abelian groups and metacyclic $p$-groups of composite order (in particular, for finite $p$-groups with $|\Omega(G)|=p^2$, $p>2$), $G$ and ${\cal K}_1(G)$ are uniquely determined by each other, up to isomorphism.

Recall that the exponent of a finite group $G$, denoted by $\exp G$, is the least common multiple of the orders of all elements of $G$.

\begin{Lemma}\label{exp}
\rm  Let $G$ be a finite group of composite order and $K$ a non-trivial subgroup of $G$. If $\varphi$ is an isomorphism from ${\cal K}_1(G)$ to ${\cal K}_1(\widetilde{G})$ for some group $\widetilde{G}$, then $|K|=|\varphi(K)|$ and $\exp K=\exp \varphi(K)$.
\end{Lemma}

\begin{Proof}
\rm Firstly, we can see from Lemma \ref{order} that $|G|=|\widetilde{G}|=|\varphi(G)|$, since $G$ is the zero of ${\cal K}_1(G)$. Also, we have that $|K|=|\varphi(K)|$, since $|K|=|G|/|G:K|=|G|/|R_{\scriptscriptstyle K}|$. Let $g$ be a non-trivial element of $K$, then $\langle g\rangle$ is a cyclic subgroup of $G$. It follows from Lemma \ref{cyclic} that $\varphi(\langle g\rangle)$ is cyclic. Suppose that $\varphi(\langle g\rangle)=\langle \widetilde{g}\rangle$. Then we have $|\widetilde{g}|=|g|$, since $|\langle \widetilde{g}\rangle|=|\langle g\rangle|$. It follows that $\exp K=\exp \varphi(K)$, as required.
\end{Proof}

\begin{Theorem}\label{isoab}
\rm  Let $G$ be a finite abelian group of composite order. Then for any group $\widetilde{G}$,
${\cal K}_1(G)\cong {\cal K}_1(\widetilde{G})$ if and only if $G\cong \widetilde{G}$.
\end{Theorem}

\begin{Proof}
\rm It is evident that $G\cong \widetilde{G}$ implies that  ${\cal K}_1(G)\cong {\cal K}_1(\widetilde{G})$. Conversely, suppose that ${\cal K}_1(G)\cong {\cal K}_1(\widetilde{G})$.
By Theorem \ref{isom1}, Theorem \ref{omega1} and \cite[4.2.6]{robinson}, we only need to consider the case where $G$ is a $p$-group and $|\Omega(G)|=p^2$, and so we may write $G=\langle a\rangle\times \langle b\rangle$, where $|a|=p^m$, $|b|=p^n$ and $m\geq n\geq 1$.

We can see from Lemma \ref{exp} that $|\widetilde{G}|=|G|$ and $\exp \widetilde{G}=\exp G=p^m$. Hence $\widetilde{G}$ is a $p$-group. Let $\widetilde{g}\in \widetilde{G}$ such that $|\widetilde{g}|=p^m$. Then we may write $\widetilde{G}=\langle \widetilde{g}\rangle\times \widetilde{A}$ for some subgroup $\widetilde{A}$ of $\widetilde{G}$ by \cite[4.2.7]{robinson}. By Lemma \ref{cyclic}, we have that $$\Omega(G)=\langle H\leq G: |H|=p\rangle=\textstyle\bigwedge\{H\in E({\cal K}_1(G)): |E(H^\upharpoonright)|=1\}.$$
This implies that the image of $\Omega(G)$ under the isomorphism from ${\cal K}_1(G)$ to $ {\cal K}_1(\widetilde{G})$ is precisely $\Omega(\widetilde{G})$ and so $|\Omega(\widetilde{G})|=|\Omega(G)|=p^2$ by Lemma \ref{order} and Lemma \ref{filter}. It follows that $\widetilde{A}$ is cyclic and $|\widetilde{A}|=|\widetilde{G}|/|\langle \widetilde{g}\rangle|=|G|/p^m=p^n$. Thus $G\cong \widetilde{G}$, as required.
\end{Proof}

\begin{Lemma}\label{centre}
\rm  Let $G$ be a finite group and $\{M_1,M_2,\ldots,M_n\}$ the set of all maximal abelian subgroups of $G$. Then $Z(G)=\bigcap_{i=1}^nM_i$.
\end{Lemma}

\begin{Proof}
\rm It is easy to see that $Z(G)\leq M_i$ for each $i\in\underline{n}$, since $Z(G)M_i$ is an abelian subgroup of $G$. Thus $Z(G)\leq \bigcap_{i=1}^nM_i$. On the other hand, for each $y\in G$, there exists $i\in\underline n$ such that $\langle y\rangle\leq M_i$, since $\langle y\rangle$ is an abelian subgroup of $G$. Thus for each $x\in\bigcap_{i=1}^nM_i$, $xy=yx$, since $x,y\in M_i$. This shows that $Z(G)\geqslant \bigcap_{i=1}^nM_i$. We have proved that $Z(G)=\bigcap_{i=1}^nM_i$, as required.
\end{Proof}

\begin{Lemma}\label{isoc}
\rm  Let $G$ be a finite group of composite order.  Then ${\cal K}_1(G)\cong {\cal K}_1(\widetilde{G})$ implies that $Z(G)\cong Z(\widetilde{G})$ for any group $\widetilde{G}$.
\end{Lemma}

\begin{Proof}
Let $\varphi: {\cal K}_1(G)\to {\cal K}_1(\widetilde{G})$ be a semigroup isomorphism. We first show that $H$ is an abelian subgroup of $G$ if and only if $\varphi(H)$ is an abelian subgroup of $\widetilde{G}$.  Suppose that $H$ is abelian subgroup of $G$. If $|H|$ a prime, then $|E(H^\upharpoonright)|=1$, and so $|E(\varphi(H)^\upharpoonright)|=1$. This implies that $|\varphi(H)|$ is also a prime and so $\varphi(H)$ is an abelian subgroup of $\widetilde{G}$. Otherwise, suppose that $|H|$ is composite. Then
$${\cal K}_1(H)=H^\upharpoonright\cong \varphi(H)^\upharpoonright={\cal K}_1(\varphi(H)),$$
since $\varphi$ preserves the natural partial order. It follows from Theorem \ref{isoab} that $\varphi(H)\cong H$, and so $\varphi(H)$ is an abelian subgroup of $\widetilde{G}$. Dually, we can show that $H$ is an abelian subgroup of $G$ if $\varphi(H)$ is an abelian subgroup of $\widetilde G$.

Let $M=\{M_1,M_2,\ldots, M_s\}$ be the set of all maximal abelian subgroups of $G$. Then
$$\widetilde M=\{\varphi(M_1),\varphi(M_2),\ldots, \varphi(M_s)\}$$
is the set of all maximal abelian subgroups of $\widetilde G$, since $\varphi$ preserves the natural partial order.
If $\bigvee M$ exists, then we can see from Lemma \ref{centre} that $Z(G)=\bigvee M$, and so
$$\varphi(Z(G))=\varphi(\textstyle\bigvee M)=\textstyle\bigvee(\varphi (M))=Z(\widetilde G).$$
This implies that
$${\cal K}_1(Z(G))=Z(G)^\upharpoonright\cong \varphi(Z(G))^\upharpoonright={\cal K}_1(\varphi(Z(G)))={\cal K}_1(Z(\widetilde G)).$$
Also, we can see from Lemma \ref{exp} that $|Z(G)|=|Z(\widetilde G)|$.
Thus, it follows from Theorem \ref{isoab} that $Z(\widetilde G)=\varphi(Z(G))\cong Z(G)$.

Finally, if $\bigvee M$ does not exist, then neither does $\bigvee \widetilde M$. In this case, both $Z(G)$ and $Z(\widetilde G)$ are trivial, and so $Z(G)\cong Z(\widetilde G)$, as required.
\end{Proof}

A group $G$ is called metacyclic if $G/N$ is cyclic for some normal cyclic subgroup $N$ of $G$. Based on King's characterization of metacyclic groups in \cite{king}, we can show that

\begin{Theorem}\label{metacy}
\rm  Let $G$ be a finite metacyclic $p$-group of composite order. Then for any group $\widetilde{G}$, ${\cal K}_1(G)\cong {\cal K}_1(\widetilde{G})$ if and only if $G\cong \widetilde{G}$.
\end{Theorem}

\begin{Proof}
\rm It is evident that $G\cong \widetilde{G}$ implies that ${\cal K}_1(G)\cong {\cal K}_1(\widetilde{G})$. Conversely, suppose that ${\cal K}_1(G)\cong {\cal K}_1(\widetilde{G})$.
We can see from Lemma \ref{kgjiben}, Lemma \ref{cyclic} and Lemma \ref{isoc} that $\widetilde{G}$ is  metacyclic and $Z(G)\cong Z(\widetilde{G})$. Also, we have that $G/G'\cong \widetilde{G}/\widetilde{G}'$, since $G/G'$ is the greatest abelian factor group of $G$, i.e., largest abelian ${\cal H}$-class which contains a central idempotent in ${\cal K}_1(G)$. Similarly, $G$ and $\widetilde{G}$ have isomorphic greatest dihedral factor groups. Further, we can see from Lemma \ref{minmax} and Lemma \ref{exp} that the minimum of the orders of the non-Frattini elements of $G$ and $\widetilde{G}$ are equal.
Thus, $G\cong \widetilde{G}$ by \cite[Theorem 5.4]{king}.
\end{Proof}

Let $G$ be a finite $p$-group with $|\Omega(G)|=p^2$, $p>2$. Noticing that $\Omega(G)$ contains all minimal subgroups of $G$, we have that $\Omega(G)$ is the unique subgroup of type $(p,p)$ in $G$. Then by \cite[Introduction, Lemma 5]{berk}, there exists  a minimal subgroup $Z$ of $G$ such that $C_G(Z)=G$. It follows from \cite[Proposition 13.26]{berk} that $G$ is metacyclic. As a consequence of Theorem \ref{metacy}, we have

\begin{Corollary}
\rm  Let $G$ be a finite $p$-group with $|\Omega(G)|=p^2$, $p>2$. Then for any group $\widetilde{G}$, ${\cal K}_1(G)\cong {\cal K}_1(\widetilde{G})$ if and only if $G\cong \widetilde{G}$.
\end{Corollary}

\noindent{\bf Remark}
\rm We have shown that $G$, ${\cal K}_1(G)$ and $C({\cal K}_1(G))$ are uniquely determined by each other, up to isomorphism, if $G$ is one of the groups appearing in Theorem \ref{isom1}, Theorem \ref{isoab} and Theorem \ref{metacy}.

\end{document}